\documentclass[12pt]{amsart}
\usepackage{graphicx}
\numberwithin{equation}{section} 

\newtheorem{thm}{Theorem}[section]
\newtheorem{cor}[thm]{Corollary}
\newtheorem{lem}[thm]{Lemma}

\newtheorem{prob}[thm]{Problem}
\theoremstyle{definition}

\theoremstyle{remark}

\newcommand{\bea}{\begin{eqnarray}}
\newcommand{\eea}{\end{eqnarray}}
\newcommand{\ba}{\begin{array}}
\newcommand{\ea}{\end{array}}
\newcommand{\bc}{\begin{center}}
\newcommand{\ec}{\end{center}}
\newcommand{\be}{\begin{equation}}
\newcommand{\ee}{\end{equation}}

\def\bn{{\mathbb N}}
\def\s{\sigma}
\def\e{{\bf 1}\!\!{\rm I}}
\def\l{\lambda}

\def\t{\tau}
\def\xb{{\mathbf{x}}}
\def\cf{{\mathcal F}}

\def\br{\mathbb{R}}
\def\a{\alpha}

\def\m{\mu}

\begin{document}

\title[ON DOMINANT CONTRACTIONS]
{A NOTE ON DOMINANT CONTRACTIONS OF JORDAN ALGEBRAS}

\author{Farrukh Mukhamedov}
\address{Farrukh Mukhamedov\\
 Department of Computational \& Theoretical Sciences\\
Faculty of Science, International Islamic University Malaysia\\
P.O. Box, 141, 25710, Kuantan\\
Pahang, Malaysia} \email{{\tt far75m@yandex.ru} {\tt
farrukh\_m@iiu.edu.my}}

\author{Seyit Temir}
\address{Seyit Temir\\
Department of Mathematics\\
Arts and Science Faculty\\
Harran University, Sanliurfa, 63200, Turkey} \email{{\tt
seyittemir38@yahoo.com}}

\author{Hasan Akin}
\address{Hasan Akin\\
Department of Mathematics\\
Arts and Science Faculty\\
Harran University, Sanliurfa, 63200, Turkey} \email{{\tt
akinhasan@harran.edu.tr}}

\begin{abstract} In the paper we consider
two positive contractions $T,S:L_{1}(A,\tau)\longrightarrow
L_{1}(A,\tau)$ such that $T\leq S$, here $(A,\t)$ is a semi-finite
$JBW$-algebra. If there is an $n_{0}\in\mathbb{N}$ such that
$\|S^{n_{0}}-T^{n_{0}}\|<1$. Then we prove that $\|S^{n}-T^{n}\|<1$
holds for every $n\geq n_{0}.$ \vskip 0.3cm \noindent


\end{abstract}

\maketitle

\section{Introduction}

Let $(X,\cf,\m)$ be a measure space with a positive $\s$-additive
measure $\m$ and let $L_1(X,\cf,\m)$ be the usual associated real
$L_1$-space. A linear operator $T:L_1(X,\cf,\m)\to L_1(X,\cf,\m)$ is
called {\it positive contraction} if  $Tf\geq 0$ whenever $f\geq 0$
and $\|T\|\leq 1$. In \cite{Z} the following theorem was proved.

\begin{thm}\label{Z1} Let $T,S:L_1(X,{\cf},\m)\to
L_1(X,{\cf},\m)$ be two positive contractions such that $T\leq S$.
If $\|S-T\|<1$ then $\|S^n-T^n\|<1$ for all $n\in\bn$
\end{thm}

Using this formulated result Zaharapol proved so called "zero-two"
law for a positive contraction of $L^1$-space. Note that the
"zero-two" law firstly appeared in \cite{OS}.

But one can ask: what would be for the given contractions the
following equality holds $\|S-T\|=1$. In this case we could not
apply the formulated theorem. Therefore, there are two options:
\begin{enumerate}
\item[(i)] one has $\|S^n-T^n\|=1$ for all $n\in\bn$;

\item[(ii)] there is an $n_0\in\bn$ such $\|S^{n_0}-T^{n_0}\|<1$.
\end{enumerate}

So, concerning (ii) we can formulate the following

\begin{prob}\label{pr} Let $S,T$ be as above in Theorem \ref{Z1}. If
there is an $n_0\in\bn$ such $\|S^{n_0}-T^{n_0}\|<1$, then can we
state $\|S^n-T^n\|<1$ for every $n\geq n_0$?
\end{prob}

By denoting $\tilde S=S^{n_0}$,$\tilde T=T^{n_0}$ as a direct
consequence of Theorem \ref{Z1} we get that $\|(\tilde S)^n-(\tilde
T)^n\|<1$ for every $n\in\bn$ under the statement of problem. This
means that $\|S^{nn_0}-T^{nn_0}\|<1$ for every $n\in\bn$. But this
is not an answer to the question.

The aim of this paper is to give an affirmative answer to the
formulated problem  for positive $L_1$-contractions of
$JBW$-algebras (in Remark 3.3 we point out that the result can be proved for
any partially ordered Banach spaces in which the norm has the additivity property).
Such a result will include as a particular case of
the Zaharopl's result. Further, we shall show that indeed that our
result is an extension of Theorem \ref{Z1}. Namely, we provide an
example of two positive contractions  for which the condition of
Theorem \ref{Z1} is not satisfied, but the statement of the problem
holds. Note that Jordan Banach algebras \cite{HS},\cite{S} are a
non-associative real analogue of von Neumann algebras. The existence
of exceptional $JBW$-algebras does not allow one to use the ideas
and methods from von Neumann algebras. To ergodic type theorems for
Jordan algebras were devoted a lot papers (see for example,
\cite{A1},\cite{A2},\cite{GKS},\cite{MTA} e.c.t.). It would be worth to mention that
a book \cite{E} is devoted to asymptotic analysis of $L_1$-contractions on commutative
and non-commutative setting. The motivation of
these investigations arose in quantum statistical mechanics and
quantum field theory (see \cite{BR}, \cite{R}). We hope that our
result will serve to prove  the "zero-two" law in a non-associative
or non-commutative framework, since nowadays such activities has
been reviewed by many authors (see for example, \cite{J1}) motivated
by various physical reasons.

\section{Preliminaries}

In this section recall some well known facts concerning Jordan
algebras.

Let $A$ be a linear space $A$ over the reals $\br$. A pair
$(A,\circ)$, where $\circ$ is a binary operation (i.e.
multiplication), is called {\it Jordan algebra} if the following
conditions are satisfied:
\begin{enumerate}
\item[(i)] $a\circ (b+c)=a\circ b+a\circ c$; $(b+c)\circ a=b\circ a+c\circ a$ for any $a,b,c\in
A$;

\item[(ii)] $\l(a\circ b)=(\l a)\circ b=a\circ(\l b)$ for any
$\l\in\br$, $a,b\in A$;

\item[(iii)] $a\circ b=b\circ a$ for any $a,b\in A$;

\item[(iv)] $a^2\circ (b\circ a)=(a^2\circ b)\circ a$ for any
$a,b\in A$.
\end{enumerate}

Let $A$ be a Jordan algebra with unity $\e$ and at the same time be
a Banach space over the reals. If a norm on $A$ respects
multiplication so that $\|a^{2}\|=\|a\|^{2}$ and $\|a^{2}\| \leq
\|a^{2}+b^{2}\|$ for all $a, \ b \in A $, then $A$ is called a {\it
$JB$-algebra} (see \cite{A3},\cite{ARU},\cite{HS}). Note that in
each $JB$-algebra $A$ the set $A^{+}= \{ a^{2}: a \in A \}$ is
regular convex cone and defines in $A$ a partial ordering compatible
with the algebraic operations. A $JB$-algebra $A$ is called a {\it
$JBW$-algebra} if there exists a Banach space $N$, which is said to
be pre-dual to $A$, such that $A$ is isometrically isomorphic to the
space $N^{*}$ of continuous linear functionals on $N$. So, on the
$JBW$-algebra $A$ one can introduce the $\sigma(A,N)$-weak topology.
It is known that the pre-dual space $N$ of a $JBW$-algebra $A$ can
be identified with the space of all $\sigma(A,N)$-weak continuous
linear functionals $A_*$ on $A$.

Recall that a  {\it trace} on a $JBW$-algebra is a map
$\tau:A^{+}\rightarrow [0,\ \infty]$ such that
\begin{enumerate}
   \item[(1)] $\tau(a+\lambda b)=\tau(a)+\lambda \tau(b)$ for all $a,\ b
   \in A^{+}$ and $\lambda\in \mathbb{R}_{+}$, provided that $0\cdot(\infty)=0$,
   \item[(2)] $\tau(U_{s}a)=\tau(a)$ for all $a\in A^{+}$ and $s\in
   A$, $s^{2}=\e$, where
$U_{s}x=2s\circ(s\circ x)-s^{2}\circ x$.
\end{enumerate}

A trace $\tau$ is said to be {\it faithful} if $\tau(a)>0$ for all
$a \in A^{+}$, $a\neq 0$; it is normal if for each increasing net
$x_{\alpha}$ in $A^{+}$ that is bounded above one has $\tau(\sup
x_{\alpha})=\sup\tau(x_{\alpha})$; it is {\it semi-finite} if there
exists a net $\{b_{\alpha}\}\subset A^{+}$ increasing to $\e$ such
that $\tau(b_{\alpha})< \infty$ for all $\alpha$, and it is {\it
finite} if $\tau(\e)< \infty$.

Throughout the paper we will consider a $JBW$-algebra $A$ with a
faithful semi-finite normal trace $\tau$. Therefore, we omit this
condition from the formulation of theorems.

Given $1\leq p<\infty$, let $A_p=\{x\in A: \tau(|x|^p)<\infty\}$,
here $|x|$ denotes the modules of an element $x$. Define the map
$\|\cdot\|_{p}:A\rightarrow [0, \ \infty)$  by the formula
$\|\cdot\|_{p}=(\tau(|a|^p))^{1/p}$. Then a pair $(A_p,\|\cdot\|_p)$
is a normed space (see \cite{A3}). Its completion in the norm
$\|\cdot\|_{p}$ will be denoted by $L_{p}(A,\t)$. As usual, we set
$L_\infty(A,\t)=A$ equipped with the norm of $A$. It is shown
\cite{A3} that the spaces $L_{1}(A,\tau)$ and $A_*$ are
isometrically isomorphic, therefore they can be indentified. Further
we will use this fact without noting.

\begin{thm}\cite{A3} The space $L_{p}(A,\tau)$, $p\geq 1$
coincides with the set
$$
L_{p}=\bigg\{ x=\int^{\infty}_{-\infty}\lambda de_{\lambda}: \
\int^{\infty}_{-\infty}|\lambda|^p d \tau (e_{\lambda})< \infty
\bigg\}.
$$
Moreover,
$$
\|x\|_{p}=\bigg(\int^{\infty}_{-\infty}|\lambda |^pd \tau
(e_{\lambda})\bigg)^{1/p}.
$$
\end{thm}

For  more information about Jordan algebras we refer a reader to
\cite{A3},\cite{ARU},\cite{HS}.

In the sequel we shall work with mappings of $L_1$-space. Therefore,
recall that a linear bounded operator $T:L_1(A,\t)\to L_1(A,\t)$ is
{\it positive} is $Tx\geq 0$ whenever $x\geq 0$. A linear operator
$T$ is said to be a {\it contraction} if $\|T\|\leq 1$. Here $\|T\|$
is defined as usual, i.e. $\|T\|=\sup\{\|Tx\|_1:\ \|x\|_1=1\}$.

\section{Main results}

In this section we are going to prove a main result of the paper.
But before do it we need some auxiliary lemmas.

\begin{lem}\label{3.1}  Let $T: L_{1}(A,\tau)\rightarrow
L_{1}(A,\tau)$ be a positive operator. Then
$$\|T\|=\sup_{\|x\|=1}\|Tx\|=\sup_{\|x\|=1,x\geq
0}\|Tx\|.$$ \end{lem}

\begin{proof} Denote $\alpha=\sup\limits_{\|x\|=1,x\geq
0}\|Tx\|$. It is clear that $\alpha\leq\|T\|$. Let $x \in
L_{1}(A,\t)$, $\|x\|=1$, then $x=x^{+}-x^{-}$,
$\|x\|=\|x^{+}\|+\|x^{-}\|$; we have \bea \|Tx\| & = &
\|Tx^{+}-Tx^{-}\|\nonumber\\
& = & \left\|\|x^{+}\|T\bigg(\frac{x^{+}}{\|x^{+}\|}\bigg)-\|x^{-}\|
T \bigg(\frac{x^{-}}{\|x^{-}\|}\bigg)\right\|\nonumber \\
&\leq
&\|x^{+}\|\left\|T\bigg(\frac{x^{+}}{\|x^{+}\|}\bigg)\right\|+\|x^{-}\|
\left\|T \bigg(\frac{x^{-}}{\|x^{-}\|}\bigg)\right\|\nonumber
\\
& \leq & \|x^{-}\|\alpha+\|x^{-}\|\a=\alpha. \nonumber \eea
Therefore $\|Tx\|\leq\alpha,$ hence $\alpha=\|T\|.$
\end{proof}

\begin{lem}\label{3.2} Let $T, S: L_{1}(A,\tau)\rightarrow
L_{1}(A,\tau)$ be two positive contraction such that $T\leq S$. Then
for every $x\in L_1(A,\t)$, $x\geq 0$ the equality holds
$$\|Sx-Tx\|=\|Sx\|-\|Tx\|.$$ \end{lem}

\begin{proof}  Let $x\in L_1(A,\t)$,$x\geq 0$, then we have \bea
\|(S-T)x\|& = &\tau
(Sx-Tx)\nonumber \\
& = & \tau(Sx)-\tau(Tx)\nonumber \\
& = & \| Sx\|-\| Tx\|. \nonumber
\eea
\end{proof}

Now we are ready to formulate the result.

\begin{thm}\label{ZN} Let
$T,S:L_{1}(A,\tau)\longrightarrow L_{1}(A,\tau)$ be two positive
contractions such that $T\leq S$. If there is an
$n_{0}\in\mathbb{N}$ such that $\|S^{n_{0}}-T^{n_{0}}\|<1$. Then
$\|S^{n}-T^{n}\|<1$ for every $n\geq n_{0}.$
\end{thm}

\begin{proof} Let us assume that $\|S^{n}-T^{n}\|=1$ for some $n>n_{0}.$
Therefore denote
\begin{equation*}
m=\min\{n\in\mathbb{N}:\|S^{n_{0}+n}-T^{n_{0}+n}\|=1\}.
\end{equation*}
It is clear that $m\geq 1$. The inequality $T\leq S$ implies that
$S^{n_{0}+m}-T^{n_{0}+m}$ is a positive operator. Then according to
Lemma \ref{3.1} there exists a sequence $\{x_{n}\}\in L_{1}(A,\t)$ such
that $x_{n}\geq 0$, $\|x_{n}\|=1, \forall n\in\mathbb{N}$ and
\begin{eqnarray}\label{eq1}
\lim\limits_{n\to\infty}\|(S^{n_{0}+m}-T^{n_{0}+m})x_{n}\|=1.
\end{eqnarray}

Positivity of $S^{n_{0}+m}-T^{n_{0}+m}$ and $x_{n}\geq 0$ together
with  Lemma \ref{3.2} imply that
\begin{eqnarray}\label{eq2}
\|(S^{n_{0}+m}-T^{n_{0}+m})x_{n}\|=\|S^{n_{0}+m}x_{n}\|-\|T^{n_{0}+m}x_{n}\|
\end{eqnarray}
for every $n\in\bn$. It then follows from \eqref{eq1},\eqref{eq2}
that

\begin{eqnarray}\label{eq3}
&&\lim\limits_{n\to\infty}\|S^{n_{0}+m}x_{n}\|=1,\\
\label{eq4} &&\lim\limits_{n\to\infty}\|T^{n_{0}+m}x_{n}\|=0.
\end{eqnarray}

The contractivity of $T$ and $S$ implies that
$\|T^{n_{0}+m-1}x_{n}\|\leq 1$, $\|T^{m}x_{n}\|\leq 1$ and
$\|S^{n_{0}}T^{m}x_{n}\|\leq 1$ for every $n\in\mathbb{N}$.
Therefore we may choose a subsequence $\{x_{n_{k}}\}$ of $\{x_{n}\}$
such that the sequences $\{\|T^{n_{0}+m-1}x_{n_{k}}\|\}$,
$\{\|T^{m}x_{n_{k}}\|\}$, $\{\|S^{n_{0}}T^{m}x_{n_{k}}\|\}$
converge. Put $y_{k}=x_{n_{k}}, k\in\mathbb{N}$ and
\begin{eqnarray}\label{eq5}
&&\alpha=\lim\limits_{k\to\infty}\|T^{n_{0}+m-1}y_{k}\|,\\
\label{eq6}
&&\beta=\lim\limits_{k\to\infty}\|S^{n_{0}}T^{m}y_{k}\|,\\
\label{eq7} &&\gamma=\lim\limits_{k\to\infty}\|T^{m}y_{k}\|.
\end{eqnarray}

From the inequalities $\|S^{n_{0}+m}y_{k}\|\leq
\|S^{n_{0}+m-1}y_{k}\|$, $\|S^{n_{0}+m}y_{k}\|\leq \|S^{m}y_{k}\|$
together with  \eqref{eq3} one gets
\begin{eqnarray}\label{eq8}
&&\lim\limits_{k\to\infty}\|S^{n_{0}+m-1}y_{k}\|=1,\\
\label{eq9} &&\lim\limits_{k\to\infty}\|S^{m}y_{k}\|=1.
\end{eqnarray}

The inequality $\|S^{n_{0}+m-1}x_{n}-T^{n_{0}+m-1}x_n\|<1$ with
\eqref{eq8} implies that $\alpha>0$. Hence  we may choose a
subsequence $\{z_{k}\}$ of $\{y_{k}\}$ such that
$T^{n_{0}+m-1}z_{k}\neq 0$, $k\in\mathbb{N}$.

Now from $\|T^{n_{0}+m-1}z_{k}\|\leq \|T^{m}z_{k}\|$ together with
\eqref{eq5}, \eqref{eq7} we find $\alpha\leq\gamma$, and hence
$\gamma>0.$

Using Lemma \ref{3.2} one gets
\begin{eqnarray}\label{eq10}
\|S^{n_{0}}T^{m}z_{k}\|&=&\|S^{n_{0}+m}z_{k}-(S^{n_{0}+m}z_{k}-S^{n_{0}}T^{m}z_{k})\|\nonumber\\
&=&\|S^{n_{0}+m}z_{k}\|-\|S^{n_{0}+m}z_{k}-S^{n_{0}}T^{m}z_{k}\|\nonumber\\
&\geq&
\|S^{n_{0}+m}z_{k}\|-\|S^{m}z_{k}-T^{m}z_{k}\|\nonumber\\
&=&\|S^{n_{0}+m}z_{k}\|-\|S^{m}z_{k}\|+\|T^{m}z_{k}\|
\end{eqnarray}
Due to \eqref{eq3},\eqref{eq9} we have
$$\lim\limits_{k\to\infty}\|S^{n_{0}+m}z_{k}\|-\|S^{m}z_{k}\|=0;$$
which with \eqref{eq10} implies that
$$\lim\limits_{k\to\infty}\|S^{n_{0}}T^{m}z_{k}\|\geq\lim\limits_{k\to\infty}\|T^{m}z_{k}\|,$$
therefore, $\beta\geq\gamma.$

On the other hand, by $\|S^{n_{0}}T^{m}z_{k}\|\leq\|T^{m}z_{k}\|$
one gets $\gamma\geq\beta$, hence $\gamma=\beta$.

Now  set
$$u_{k}=\frac{T^{m}z_{k}}{\|T^{m}z_{k}\|}, \ \ k\in\bn.$$
Then using the equality $\gamma=\beta$ and \eqref{eq4}  one has
\begin{eqnarray*}
&&\lim\limits_{k\to\infty}\|S^{n_{0}}u_{k}\|=
\lim\limits_{k\to\infty}\frac{\|S^{n_{0}}T^{m}z_{k}\|}{\|T^{m}z_{k}\|}=1,\\
&&\lim\limits_{k\to\infty}\|T^{n_{0}}u_{k}\|=
\lim\limits_{k\to\infty}\frac{\|T^{n_{0}+m}z_{k}\|}{\|T^{m}z_{k}\|}=0.
\end{eqnarray*}

So, owing to Lemma \ref{3.2} and positivity of
$S^{n_{0}}-T^{n_{0}}$, we get

$$\lim\limits_{k\to\infty}\|(S^{n_{0}}-T^{n_{0}})u_{k}\|=1.$$
Since $\|u_{k}\|=1, u_{k}\geq 0, \forall k\in\mathbb{N}$ from Lemma
\ref{3.1} one finds $\|S^{n_{0}}-T^{n_{0}}\|=1,$ which is a
contradiction. This completes the proof.
\end{proof}

As a corollary of the proved theorem we obtain Zaharopol's result in
a non-associative setting. Moreover, it recovers one when the
algebra is associative.

\begin{cor}\label{Z}  Let $T, S: L_{1}(A,\tau)\rightarrow
L_{1}(A,\tau)$ be two positive contractions such that $T\leq S$. If
$\| S-T \| <1,$ then  $\| S^{n}-T^{n}\| < 1$ for every $n\geq 1$.
\end{cor}

{\bf Remark 3.1.}  It should be noted the following:
\begin{enumerate}
\item[(i)] Since the dual of $L^1(A,\t)$ is $A$, then due to the
duality theory the proved Theorem \ref{ZN} holds if we replace
$L^1$-space with $JBW$-algebra $A$.

\item[(ii)] Unfortunately, Theorem \ref{ZN} is not longer true if one
replaces $L_1$-space by  an $L_p$-space, $1<p<\infty$. The
corresponding example was provided in \cite{Z}.

\item[(iii)] It would be better to note that certain ergodic properties
of dominant positive operators has been studied in \cite{EW} in a
non-commutative setting. In general, to dominant operators were
devoted a monograph \cite{K}.

\end{enumerate}

 Note that this corollary does not imply the
proved theorem. Indeed, let us consider the following

{\bf Example.} Consider $\br^2$ with a norm $\|\xb\|=|x_1|+|x_2|$,
where $\xb=(x_1,x_2)$. An order in $\br^2$ is defined as usual,
namely $\xb\geq 0$ if and only if $x_1\geq 0$, $x_2\geq 0$. Now
define mappings $T:\br^2\to\br^2$ and $S:\br^2\to\br^2$,
respectively, by
\begin{eqnarray}\label{S}
&& S(x_1,x_2)=(Ax_1+Bx_2,Cx_1+Dx_2),\\
\label{T} && T(x_1,x_2)=(\l x_2,0).
\end{eqnarray}
The positivity of $S$ and $T$ implies that $A,B,C,D,\l\geq 0$. It is
easy to check that $T\leq S$ holds if and only if $\l\leq B$.

One can see that
\begin{eqnarray}\label{S2}
 S^2(x_1,x_2)&=&\big((A^2+BC)x_1+(AB+BD)x_2,\nonumber
\\&&(AC+DC)x_1+(D^2+BC)x_2\big),\\
\label{T2}  T^2(x_1,x_2)&=&(0,0).
\end{eqnarray}

By means of Lemma \ref{3.1} let us calculate the norms of operators
$S$,$S-T$,$S^2-T^2$. Furthermore, we assume that
\begin{equation}\label{AC}
B+D\leq A+C.
\end{equation}

Then using \eqref{S},\eqref{T} we have
\begin{eqnarray}\label{S1}
\sup_{\|\xb\|=1\atop \xb\geq 0}\|S\xb\|&=&\max_{x_1+x_2=1\atop
x_1,x_2\geq 0}
\{(A+C)x_1+(B+D)x_2\}\nonumber\\
&=&\max_{0\leq x_1\leq 1}\{(A+C-B-D)x_1+B+D\}\nonumber\\
&=&A+C
\end{eqnarray}
here we have used \eqref{AC}.

Similarly, one has
\begin{eqnarray}\label{S-T}
\sup_{\|\xb\|=1\atop \xb\geq 0}\|(S-T)\xb\|&=&
\max_{0\leq x_1\leq 1}\{(A+C-B-D+\l)x_1+B+D-\l\}\nonumber\\
&=&A+C
\end{eqnarray}

Finally using \eqref{S2}, \eqref{T2} together with \eqref{AC} we
obtain
\begin{eqnarray}\label{S-T2}
\sup_{\|\xb\|=1\atop \xb\geq 0}\|(S^2-T^2)\xb\| &=&
\max_{0\leq x_1\leq 1}\bigg\{(A+D)(A+C-B-D)x_1\nonumber\\
&&+D^2+AB+BD+BC\bigg\}\nonumber\\
&=&A^2+AC+BC+DC.
\end{eqnarray}

Now from \eqref{S1},\eqref{S-T} we conclude that the equality
$A+C=1$ implies the contractivity of $S$ and $\|S-T\|=1$. The
condition \eqref{AC} yields that $T$ is a contraction.

The condition $\|S^2-T^2\|<1$ due to \eqref{S-T2} yields that
$$
A^2+AC+BC+DC<1
$$
which together with $A+C=1$ implies that $C>0$ and $B+D<1$.

Base on the finding conditions let us provide more concrete example,
i.e. $A=C=1/2$, $B=D=1/3$ and $\l=1/4$.

So, we have constructed two positive contractions $T$ and $S$ with
$S\geq T$ such that $\|S-T\|=1$,$\|S^2-T^2\|<1$. This shows that the
condition of Corollary \ref{Z} is not satisfied, but due Theorem
\ref{ZN} we have $\|S^n-T^n\|<1$ for all $n\geq 2$. Therefore the
proved Theorem \ref{ZN} is an extension of the Zaharopol's result.\\

 {\bf Remark 3.2.} Let $M$ be a  von Neumann
algebra with normal faithful  semi-finite trace $\t$ (see \cite{BR}
for definitions). By $M_{sa}$ we denote the set of all self-adjoint
elements of $M$. Let $\a: M\to M$ be a positive linear operator. A
linear operator $\a$ is said to be {\it absolute contraction} if
$\t(\a(x))\leq\t(x)$ for all $x\geq 0$,$x\in M$ and $\a(\e)\leq\e$.
Let $L_1(M,\t)$ be $L_1$-space associated with $M$. Then it is known
\cite{Y} that any absolute contraction can be extended to
$L_1(M,\t)$ such that $\|\a(x)\|\leq\|x\|$ for every $x\in
L_1(M,\t)$, $x=x^*$. We also know (see \cite{A3},\cite{ARU}) that
the self-adjoint part $M_{sa}$ of $M$ is a $JBW$-algebra with
respect to multiplication $x\circ y=(xy+yx)/2$ and
$L_1(M,\t)=L_1(M_{sa},\t)+iL_1(M_{sa},\t)$. Hence, every absolute
contraction is $L_1$-contraction of the $JBW$-algebra $M_{sa}$.
Therefore, all proved theorems will be valid for any absolute
contraction of von Neumann algebras.

{\bf Remark 3.3.} Note that Theorem \ref{ZN} can be extended to any
partially ordered Banach space $X$ in which the norm should satisfy the
{\it additivity condition} on positive part $X_+$ of $X$, i.e for any
positive elements $x_1,x_2\in X_+$ one has $\|x_1-x_2\| = \|x_1\| + \|x_2\|$.
The extended theorem's proof will remain the same as the proof of Theorem \ref{ZN}.
An example of a Banach space which has the additivity condition, besides $L_1$-spaces,
is the dual $A(K)^*$ of the space $A(K)$ of continuous affine
functions on a compact convex set $K$ (see \cite{EK}).

\section*{Acknowledgement}
The first named author (F.M.) thanks TUBITAK and Harran University
for kind hospitality. He also acknowledges Research Endowment Grant
B (EDW B 0801-58) of IIUM. Finally, the authors would like to thank the
referee for useful suggestions and remarks which improved the text of the paper.

\end{document}